# The base-*e* representation of numbers and the power law

Subhash Kak
Oklahoma State University

*Abstract*. Some properties of the optimal representation of numbers are investigated. This representation, which is to the base-e, is examined for coding of integers. An approximate representation without fractions that we call WF is introduced and compared with base-2 and base-3 representations, which are next to base-e in efficiency. Since trees are analogous to number representation, we explore the relevance of the statistical optimality of the base-e system for the understanding of complex system behavior and of social networks. We show that this provides a new theoretical explanation for the nature of the power law exhibited by many open complex systems. In specific, we show that the power law distribution most often proposed for such systems has a form that is similar to that derived from the optimal base-e representation.

**Introduction**
Novel representation schemes for numerals are important in communication engineering and computer science applications (Kak, 2016), especially when one wishes to represent data as efficiently as possible (McEliece, 2002). They are also important in understanding complex systems as manifested in behaviors such as the first digit phenomenon (Hill, 1998; Hayes, 2001), power laws and series models in social networks (Newman, 2005; Kak, 2017; Kak, 2018a), protein-protein interaction networks (De Las Rivas and Fontanillo, 2010), transportation networks (Gibbons, 1985), city and firm sizes (Axtell, 2001; Gabaix, 1999; Cristelli et al, 2012) and financial networks (Stiglitz, 2016).

We begin with number representation. Let's consider coding of numbers to an arbitrary base r. In general we should be able to write the expansion for a number *x* in the form:

$$x = \sum_{n=-\infty}^{\infty} a_n r^n \qquad (1)$$

where $a_n$ are integers so that $a_n$ is zero for a sufficiently large *n*, and $0 \leq a_n < r$. In general, this expression applies to any kind of *r*, even non-integer (Eggan and Vanden Eynden, 1966) and irrational (Bergman, 1957; Rousseau, 1995).

**Tree structure**
The polynomial representation of equation (1) may be viewed as a tree. Thus, for example, the tree for the representation for the number 43 to base 3 is derived by successive division as follows:

$$43 = 14 \times 3 + 1$$
$$14 = 4 \times 3 + 2$$
$$4 = 1 \times 3 + 1$$
$$1 = 0 \times 3 + 1$$

Therefore 43 (base 3) = $1 \times 3^3 + 1 \times 3^2 + 2 \times 3^1 + 1 = 1\ 1\ 2\ 1$, which is obtained by taking the remainders in the reverse order. This tree is shown in Figure 1, where the nodes are read from left to right.

Note that the number of branches at each node is 3, but only the named branch is shown. If one were to gather all these branches and fold them on top of each other correctly, we will have a single layered structure with the overlapped number of branches for 0, 1, and 2, respectively. One expects that aggregate numbers under each of these branches are equal and so we can assign them the same probability.

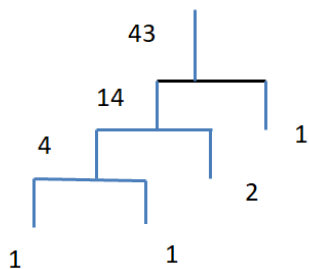

Figure 1. Tree structure of the ternary representation of 43.

Now consider the efficiency of representation based on information of the symbols. The information *I(x)* associated with the symbol *x* (which is a number less than the base *r*) is given by:

$$I(x) = -\log P(x) \qquad (2)$$

where $P(x)$ is the probability of occurrence of the symbol *x*. Equation (1) gives the information in *bits* if it uses a logarithm to base 2, and in *nats* if base e is used.

Clearly, the capacity of a base representation to carry information goes up as the size of the base increases. But the increase in information must be squared off against the extra burden entailed by the use of the larger set (which corresponds to a greater number of branches at each node). For binary, the value is $\ln 2 = 0.693$ nats (= 1 bit), whereas for *e*, it is 1 nat (=1.443 bits); for base 3, the figure is 1.099 nats (= 1.585 bits); and for base 10, it is 2.303 nats (=3.322 bits).

The probability of the use of each of the *r* symbols may be taken to be the same and equal to $1/r$, and the information associated with each symbol is $\log r$. The efficiency of the coding scheme per symbol is

$$E(r) = \frac{\ln r}{r} \qquad (3)$$

**Theorem 1**. The optimal base for number representation is *e*.

*Proof.* To find the value of *r* for which it is a maximum we take the derivative of *E(r)* with respect to *r* and equate that to zero. This yields the condition that $\ln r = 1$, from which we conclude that the optimal base is *e*, with *E*(e)= 0.368 nats or 0.531 bits.

The fact that we define efficiency based on information content means that the optimality is true only in a probabilistic sense. This indicates that comparison of base-*e* number representation to standard base-2, base-3, or decimal representations for integers alone (which is a small subset of the real numbers), or any specific example, may not always show superior results.

Table 1 below provides E(r) for some of the values of *r* that range from 2 to 10.

Table 1. Efficiency of coding for certain bases

| b | 2 | e | 3 | 4 | 5 | 8 | 10 |
|---|---|---|---|---|---|---|---|
| E(r) nats | 0.347 | 0.368 | 0.366 | 0.347 | 0.322 | 0.260 | 0.230 |
| E(r) bits | 0.500 | 0.531 | 0.528 | 0.500 | 0.465 | 0.375 | 0.331 |

The efficiency is quite close to the maximum for *r*=3 (a value of 0.528 bits as compared to the optimal value for e which is 0.531 bits), with the next best value coming at the bases 2 and 4 (where it is 0.500 bits). The efficiency at *r*=3 is superior to that at *r*=2 by 5.6% (for details, see Hill, 1998; Hayes, 2001; Kak, 2018b). After this the values decline monotonically as shown in Figure 2.

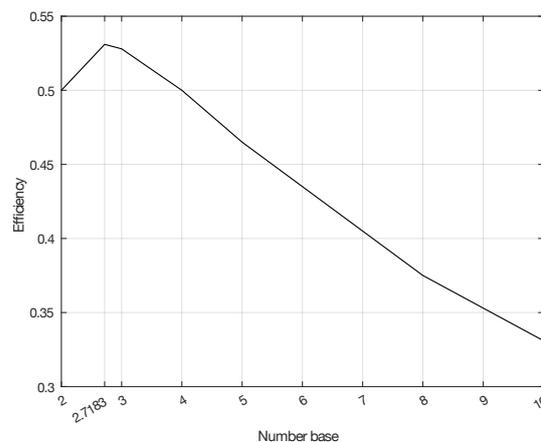

Figure 2. Efficiency of number bases, 2 through 10

In this article, we show how the base-*e* system works and compare it to the the bases of 2 and 3, which are closest to it in efficiency. Given that *e* is irrational and therefore the representations for integers in base-*e* are not going to be efficient, we further present an approximate representation that makes it easier to do a clear comparison with base-2 and -3 representations. If it is accepted that Nature chooses optimal schemes, then the base-*e* representation should show up in complex systems, such as social networks or city size, and we provide evidence that supports this hypothesis. We further show that the power law distribution often proposed for complex systems is related closely to the optimal base-*e* representation.

**The base-*e* system**

The base-*e* system represents the number using powers of *e* and coefficients, $d_i, i = \ldots 2, 1, 0, -1, -2, \ldots$, that are integers less than *e*, that is 0, 1, and 2 in the most economical manner. Clearly, owing to the nature of *e*, the representation of integers will involve decimal expansions.

$$x = d_n d_{n-1} \ldots d_2 d_1 d_0 . d_{-1} d_{-2} \ldots d_{-m} \qquad (4)$$

$$x = e^n d_n + e^{n-1} d_{n-1} + \cdots + e d_1 + d_0 + e^{-1} d_{-1} + e^{-2} d_{-2} + \cdots$$

For easy reference, the various powers of *e* are:

| n | -5 | -4 | -3 | -2 | -1 | 0 | 1 | 2 | 3 | 4 | 5 |
|---|---|---|---|---|---|---|---|---|---|---|---|
| $e^n$ | .0007 | .018 | .050 | .135 | .368 | 1 | 2.718 | 7.389 | 20.086 | 54.598 | 148.413 |

Just as 10, 100, 1000 and so on are turning points (exponentials) for base 10, the turning points for base-*e* (powers of *e*) after the numbers are rounded off are 3, 7, 20, 55, 148, 403, 1897,.. and so on.

An interesting aspect of the base-*e* system is that it maps all numbers to irrational points on the real line. In this it represents the dual to schemes with integer bases that map corresponding points to rational numbers.

In the consideration of codes for integers, the question of the degree of accuracy comes in. The superiority of the base-*e* representation is true when we consider all real numbers, but here we only wish to compare representations for integers.

In Table 2 we map integers from 1 to 20 to the base-*e* representation accurate to four "decimal" places. With the added constraint on least error for number of chosen decimal points, one can assert that a unique representation is defined.

Table 2. The base-e coding for numbers 1 through 20, accurate to four "decimal" places

| n | Base-e, to four "decimal" places | Base-e, exact value |
|---|---|---|
| 1 | 1 | 1 |
| 2 | 2 | 2 |
| 3 | 10.0200 | 2.99 |
| 4 | 11.0200 | 3.99 |
| 5 | 12.0200 | 4.99 |
| 6 | 20.1110 | 5.99 |
| 7 | 21.1110 | 6.99 |
| 8 | 22.1110 | 7.99 |
|   | 100.1120 | 7.99 |
| 9 | 101.1120 | 8.99 |
| 10 | 102.1120 | 9.99 |
| 11 | 110.2101 | 10.99 |
| 12 | 111.2101 | 11.99 |
| 13 | 112.2101 | 12.99 |
| 14 | 121.0102 | 13.99 |
| 15 | 122.0102 | 14.99 |
| 16 | 201.0110 | 15.96 |
| 17 | 202.0110 | 16.96 |
| 18 | 210.1100 | 17.99 |
| 19 | 211.1100 | 18.99 |
| 20 | 212.1100 | 19.99 |

Here's the explanation for how number 8 has two forms 22.1110 and 100.1120 in Table 2. In principle, one of these will be superior to the other in how close it is to the integer based on the number of decimal points that are chosen. Given the choice of four decimal points, 22.1110 equals 7.989, whereas 100.1120 equals 7.992; therefore, the latter representation is superior.

It is clear that base-*e* representations of numbers for integers from 1 to 20 are much less efficient than base-2 or 3 systems. Consider the $\varepsilon$-approximation of the base-*e* representation as

$$x = \sum_{i=-m}^{n} a_i\, e^i \qquad (5)$$

with the condition that

$$\left| x - \sum_{i=-m}^{n} a_i\, e^i \right| < \varepsilon$$

When $\varepsilon = 0.5$, the number is correct to the nearest integer value, and the mapping will be called the WF representation.

Table 3. The WF representation and true value

| n | Base-e, without fraction (WF) | True value |
|---|---|---|
| 1 | 1 | 1 |
| 2 | 2 | 2 |
| 3 | 10 | 2.718 |
| 4 | 11 | 3.718 |
| 5 | 12 | 4.718 |
| 6 | 21 | 6.436 |
| 7 | 22 | 7.436 |
| 8 | 101 | 8.389 |
| 9 | 102 | 9.389 |
| 10 | 110 | 10.107 |
| 11 | 111 | 11.107 |
| 12 | 112 | 12.107 |
| 13 | 120 | 12.825 |
| 14 | 121 | 13.825 |
| 15 | 122 | 14.825 |
| 16 | 201 | 15.778 |
| 17 | 202 | 16.778 |
| 18 | 211 | 18.496 |
| 19 | 212 | 19.496 |
| 20 | 220 | 20.214 |

**Base-*e* (WF) representation**

We have already mentioned that the optimality is to be understood statistically when all the points on the real line are chosen, but one would like to see it operationally for small integers. We do so by leaving out the fractions and choosing the number closest to the integer being represented as is done in quantization in signal theory (McEliece, 2002). This WF representation is mathematically:

$$x \approx \sum_{i=0}^{n} a_i \, e^i \qquad (6)$$

with the condition that

$$\left| x - \sum_{i=0}^{n} a_i \, e^i \right| < 0.5$$

This actual error between *x* and its WF representation can be either negative or positive just as is the case in quantization.

Now in Table 3, we compare the WF representation with representations for bases 2 and 3.

Table 3. Comparison of base-*e* (WF) coding against bases 2 and 3

| n | Base-e, WF | Base-2 | Base-3 |
|---|---|---|---|
| 1 | 1 | 1 | 1 |
| 2 | 2 | 10 | 2 |
| 3 | 10 | 11 | 10 |
| 4 | 11 | 100 | 11 |
| 5 | 12 | 101 | 12 |
| 6 | 21 | 110 | 20 |
| 7 | 22 | 111 | 21 |
| 8 | 101 | 1000 | 22 |
| 9 | 102 | 1001 | 100 |
| 10 | 110 | 1010 | 101 |
| 11 | 111 | 1011 | 102 |
| 12 | 112 | 1100 | 110 |
| 13 | 120 | 1101 | 111 |
| 14 | 121 | 1110 | 112 |
| 15 | 122 | 1111 | 120 |
| 16 | 201 | 10000 | 121 |
| 17 | 202 | 10001 | 122 |
| 18 | 211 | 10010 | 200 |
| 19 | 212 | 10011 | 201 |
| 20 | 220 | 10100 | 202 |

For numbers 1 through 20, we see the lengths of the representational sequences for bases *e* and 3 are identical excepting for the number 8.

The total number of characters needed is 51 for base-*e*, 74 for base-2, and 50 for base-3.

  Total cost for base-*e*: 51 × ln e = 51 nats

  Total cost for base-2: 74× ln 2 = 51.282 nats

  Total cost for base-3: 50× ln 3= 54.93 nats

Thus for the WF case for integers 1 through 20, base-*e* is superior to both the bases of 2 and 3. But as shown in an earlier paper (Kak, 2018a), the efficiency of the bases will depend on the nearness to powers of the base and there will be numbers less than the turning point of the base where that base may be superior to the others.

**Random trees and optimal branching**
We now consider trees in natural systems where in contrast to the representation of a single number as in the previous sections we look for all possible outputs associated with the system. The study of natural systems has shown a good match with heavy-tailed distributions such as Pareto distribution or the power law (Stumpf, 2012). Such systems are approximately scale-free or self-similar. Examples of these include social networks and collaboration networks (Kadushin, 2012), many kinds of computer networks, and the internet and the web graph of the World Wide Web (Newman, 2005). Preferential attachment models have been proposed as mechanisms to explain such distributions but what follows is a new approach to the problem.

In a complex system the branching lines from each node will be a random function (Newman, 2005). When the branches are folded over, let the probabilities be represented by $P(x_i)$ for the number of aggregated branches that can take values of $i$ that vary from 1 to $n$. The entropy, H(X), of the system will be:

$$H(X) = - \sum_i P(x_i) \log P(x_i) \qquad (7)$$

In the best case each of these branches will have the same probability and the mapping is most efficient. We have already seen from equation (3) that the optimal branching number is given by $e$.

For unconstrained random trees, the probability of events in ranked order is proportional to $1/n$, where n is the rank. If a counting process is uniformly distributed over the range $\{1, …, S\}$, with random values of $S$, then the sum satisfies the Newcomb-Benford Law (Hill, 1995), $P(n) = \log_r (1 + \frac{1}{n})$, where $n$ is the leading digit ($n \in \{1,2, ..., r-1\}$). When $r = e$, we have:

$$P(n) = \ln\left(1 + \frac{1}{n}\right) \qquad (8)$$

$$= \frac{1}{n+1} + \frac{1}{2(n+1)^2} + \frac{1}{3(n+1)^3} + \frac{1}{4(n+1)^4} + \cdots \qquad (9)$$

If the higher order terms are ignored, we have

$$P(n) \approx \frac{k}{n} \qquad (10)$$

This is Zipf's distribution (Belevitch, 1959; Powers, 1998; Zipf, 1949).

*Definition*. Let a counting event be a function of some aggregated property in a random tree with $r$ branches at each level.

**Theorem 2**. The number of counts $N(n+1)$ at level $n+1$ is $rN(n)$.

*Proof*. The proof is elementary because at each level there is a $r$-fold branching.

Now consider optimal branching, which by Theorem 1 is *e*-fold. The total number of branches will be limited by constraints associated with the physical system, and let's call the maximum value of the count associated with the event of interest to be *Max*.

One can read the growing values directly or read them in reverse. In the direct expansion of the tree we have starting with an initial count of *A*:

$$A \to Ae \to Ae^2 \to Ae^3 \to \cdots \to Max$$

where *Max* is the maximum count associated with the variable. Written in the reverse order, the sequence is:

$$Max \to Max\, e^{-1} \to Max\, e^{-2} \to Max\, e^{-3} \ldots \ldots \quad (11)$$

Let *k* be the variable associated with location of the values of expression (11), starting with 1 and then 2, 3, ...

Since the count expands by a factor of *e* in each branching, it will be the highest at the nodes at the bottom of the tree. When traversing the tree backwards, the count at each step will be decreased by a factor of *e*. Therefore, we can assert:

**Theorem 4.** For optimal branching in natural systems, when the aggregated counts are written in ranked order, the probability at rank *k*+1 is proportional to *ke*$^{-1}$.

This means that the scaling law is:

$$p(k) \sim e^{-(k-1)} \quad (12)$$

which is a heavy-tailed distribution.

Experiments have shown (e.g. Clauset et al., 2009) that many phenomena follow the following power law approximately for large values of *k*:

$$p(k) \sim k^{-\gamma} \quad (13)$$

where $\gamma$ is a parameter whose value is typically in the range $2 < \gamma < 3$. The main characteristic of this distribution is its heavy-tailed nature and there is considerable variety in the nature of the law.

As shown in Table 4, the values of (12) and (13) evolve in a similar fashion. For comparison, we have chosen $\gamma = 2.5$.

Table 4. Comparison of two heavy-tailed distributions of (12) and (13)

| $k$ | 1 | 2 | 3 | 4 | 5 | 6 |
|---|---|---|---|---|---|---|
| $p(k)\sim e^{-(k-1)}$ | 1 | 0.368 | 0.135 | 0.050 | 0.018 | 0.007 |
| $p(k)\sim k^{-\gamma}$ | 1 | 0.178 | 0.064 | 0.031 | 0.018 | 0.011 |

The values of the distribution in the new expression of (12) fall even faster than in (13), as it is an exponential function.

**Saturated power law**

Physical constraints will saturate the growth of the nodes down the random tree. If $\alpha$ is the saturation parameter ($\alpha < 1$), then in the growth of the tree e will be replaced by $\alpha e$. Finally, instead of formula (12) we obtain:

$$p(k)\sim(\alpha e)^{-(k-1)}, \quad \alpha < 1 \qquad (14)$$

Equation (12) may be seen as a special case of (14) when $\alpha = 1$. The $p(x)\sim e^{-(k-1)}$

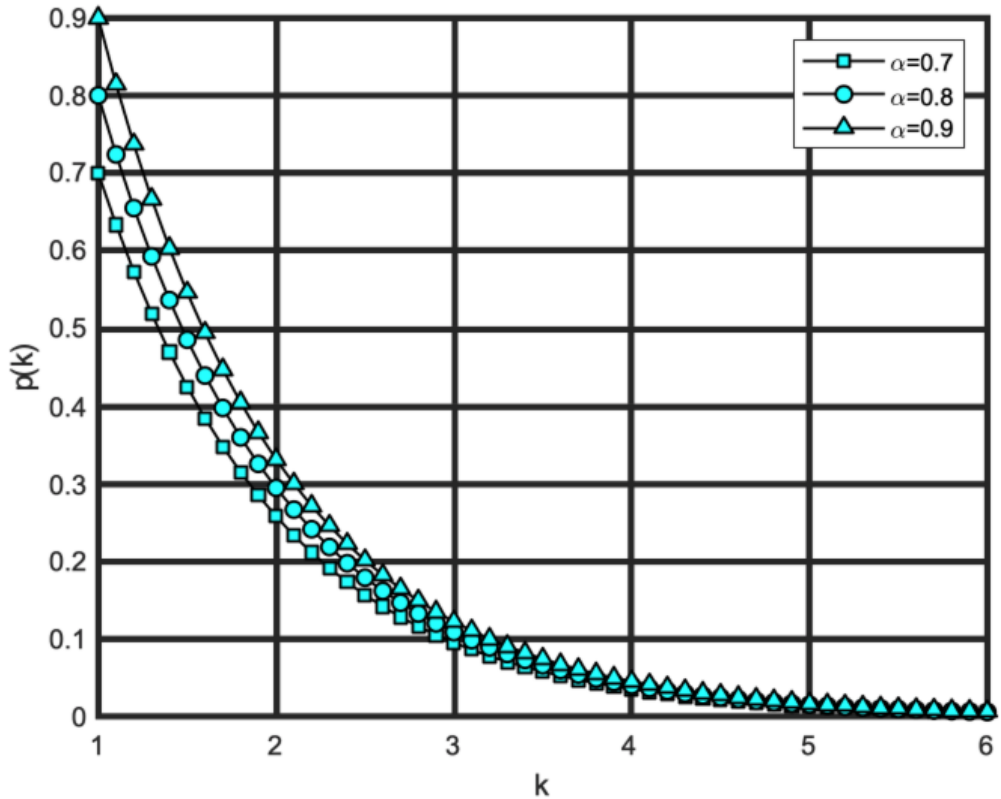

Figure 3. Saturated power law for equation (14)

The parameter in the saturated power law would vary based on the constraints associated with the natural system.

**Conclusions**

Some properties of the optimal base-*e* representation were described, and for integers an approximate representation without fractions that we call WF was introduced. Comparisons of the WF representation with those of bases-2 and 3 were made. Since trees are analogous to number representation, we explored the applicability of the base-*e* system to understanding complex system behavior. We showed that this provides a new explanation for the power law exhibited by a natural system.

The power law derived in this paper is heavy-tailed like the ones that have been widely discussed in the literature on complex systems.